\documentclass[a4paper,reqno, 12pt]{amsart}
\usepackage{amsmath,amssymb,amsfonts,amsthm, mathrsfs}
\usepackage{marginnote}
\usepackage{lmodern}
\usepackage{makecell}
\usepackage{diagbox}
\usepackage{multirow}
\usepackage{booktabs}
\usepackage{verbatim,wasysym,cite}
\usepackage{microtype}
\usepackage{color,enumitem,graphicx}
\usepackage[english]{babel}
\usepackage[symbol]{footmisc}
 \usepackage[colorlinks=true, linkcolor=blue, citecolor=red]{hyperref}
\renewcommand{\epsilon}{{\varepsilon}}
\numberwithin{equation}{section}
\newtheorem{theorem}{Theorem}[section]
\newtheorem{conjecture}{Conjecture}[section]
\newtheorem{lemma}[theorem]{Lemma}
\newtheorem{remark}[theorem]{Remark}
\newtheorem{definition}[theorem]{Definition}
\newtheorem{proposition}[theorem]{Proposition}

\begin{document}

\title[Defocusing   high dimensional   NLS]{Global well-posedness and scattering    for the defocusing  energy supercritical      NLS  in high dimensions}
	
 \author[Xuan Liu]{Xuan Liu}
\address{ School of Mathematics, Hangzhou Normal   University, Hangzhou 311121, China}
\email{liuxuan95@hznu.edu.cn }

\begin{abstract}
 We consider the defocusing energy-supercritical nonlinear Schr\"odinger equation 
  $i\partial_{t}u+\Delta u=|u|^p u$    in   dimensions  $d\ge5$. Killip-Visan [Comm. Partial Differential Equations,  2010] and Li-Li [Siam J. Math. Anal., 2022] proved that for  $s_c:=\frac{d}{2}-\frac{2}{p}>1$, any solution that remains bounded in the critical Sobolev space   $\dot H_x^{s_c}(\mathbb{R} ^d)$  must be global and scatter. In dimensions \(d \ge 8\),  their results required either that \(p\) be even or that \(s_c < \frac{d+2-\sqrt{(d-2)^2-16}}{4}\). In this paper, we improve the upper bound on \(s_c\) to \(s_c<1+p\) by establishing some new nonlinear estimates. This allows us to cover all cases in which \(p\) lies in the local existence range.
 
\vspace{0.3cm}


\end{abstract}
	
\maketitle
 
\medskip

\section{Introduction}
\subsection{Main results}

We  consider  the Cauchy probem for the defocusing    energy supercritical     nonlinear  Schr\"odinger equation    in dimension  $d\ge5$,  
\begin{equation}
\begin{cases}
i\partial_{t}u+\Delta  u=|u|^{p }u,\qquad (t,x)\in \mathbb{R} \times \mathbb{R}^d \\
u(0,x)=u_0(x) ,
\end{cases}\label{NLS}
\end{equation}
where  $p>\frac{4}{d-2}$, and    $u:\mathbb{R} _t\times \mathbb{R} _x\mapsto \mathbb{C}$  is a complex-valued function.

Equation \eqref{NLS} possesses several symmetries and invariances, among which the most important one is the scaling symmetry. That is, if $u$ solves \eqref{NLS} with initial data $u_0$, then 
  \[u_\lambda(t,x) = \lambda^{2/p } u(\lambda^2 t, \lambda x), \qquad \lambda>0  \]
  solves \eqref{NLS} with initial data $u_{0,\lambda}(x) = \lambda^{2/p } u_0(\lambda x)$ and $\|u_{0,\lambda}\|_{\dot{H}^s(\mathbb{R} ^d)} =\lambda ^{s+\frac{2}{p}-\frac{d}{2}} \|u_0\|_{\dot{H}^s(\mathbb{R} ^d)}$. The problem becomes scale-invariant when the data belongs to the homogeneous Sobolev space $\dot{H}^{d/2 - 2/\alpha }(\mathbb{R}^d)$.  This regularity index is called the \emph{critical} index  with respect to scaling and is denoted by 
 \begin{equation}
 	s_c := \frac{d}{2} - \frac{2}{p}. \label{Es}
 \end{equation}  
  
  Throughout this paper, we will restrict our consideration to the class of strong solutions. The precise notion is defined as follows.

  \begin{definition}[Strong solution]\label{Defsolution}
  	A function $u: I \times \mathbb{R}^d\to \mathbb{C}$ is a solution to \eqref{NLS} if for any compact $J \subset I$, $u \in C_t^0 \dot H_x^{s_c}(J \times \mathbb{R}^d) \cap  L_{t,x}^{\frac{d p}{2}}(J \times \mathbb{R}^d)$, and satisfies the Duhamel formula for all $t, t_0 \in I$:
  	\begin{equation}
  		u(t) = e^{i(t - t_0)\Delta} u(t_0) - i \int_{t_0}^t e^{i(t-\tau)\Delta}  (|u|^{p}u)(\tau)\, d\tau. \notag
  	\end{equation}
  	We refer to $I$ as the \textit{lifespan} of $u$. If $u$ cannot be extended to any strictly larger interval, it is called a \textit{maximal-lifespan solution}. If $I = \mathbb{R}$, $u$ is called a \textit{global solution}.
  \end{definition}

 It is well-known that, for initial data in  $\dot H^{s_c}(\mathbb{R} ^d)$, local well-posedness of (\ref{NLS}) (under the smoothness assumption $s_c<1+p $ or  $p$ is even)   can be established via Strichartz estimates and the contraction mapping principle, with a time of existence depends on the profile of the initial data.   Moreover, if the $\dot H^{s_c}(\mathbb{R} ^d)$ norm of  the initial data is sufficiently small,  then the corresponding solution exists globally in time and scatters   as $t \to \pm \infty$ in the sense that  there exist $u_{\pm} \in \dot{H}^{s_c}(\mathbb{R}^d)$ such that
 \begin{equation}\label{1.4}
 	\lim_{t \to \pm \infty} \left\| u(t) - e^{it\Delta} u_{\pm} \right\|_{\dot{H}^{s_c}_x(\mathbb{R}^d)} = 0.
 \end{equation}
  See \cite{CazenaveWeissler1990NA} for the details. 
  
We are interested in the large data global well-posedness and scattering theory for the defocusing NLS \eqref{NLS}.    It is conjectured in \cite{KenigMerle2010}  that, assuming some \textit{a priori} control of a critical norm, global well-posedness and scattering hold for any $s_c > 0$ and in any spatial dimension:
  
\begin{conjecture}\label{CNLS0}
 Let $d \geq 1$ and $p \geq \frac{4}{d}$. Assume $u: I \times \mathbb{R}^d \rightarrow \mathbb{C}$ is a maximal-lifespan solution to (\ref{NLS}) such that 
 \begin{equation}\label{E7301}
  u \in L_t^\infty \dot{H}_x^{s_c}(I \times \mathbb{R}^d),
  \end{equation}
  then $u$ is global and scatters as $t \to \pm \infty$.
  \end{conjecture}
  
 For the energy-critical ($s_c=1$) and mass-critical ($s_c=0$) cases, the conservation of mass 
 	\[ M(u(t)) := \int_{\mathbb{R}^d} |u(t,x)|^2 \, dx = M(u_0), \]
 and the energy,
 \[E(u(t)) := \int_{\mathbb{R}^d} \left(\frac{1}{2} |\nabla u(t,x)|^2 + \frac{1}{p+2} |u(t,x)|^{p +2} \right)\, dx = E(u_0) \]
   makes it unnecessary to assume (\ref{E7301}). In fact, the global well-posedness and scattering problem for large initial data in these cases had already been fully resolved. See \cite{Bourgain1999,Colliander2008,RyckmanVisan2007, Visan2007, Visan2012} for the energy-critical case, and \cite{Dodson2012,Dodson2015,Dodson2016a,Dodson2016b,TaoVisanZhang2007,KillipTaoVisan2009,KillipVisanZhang2008} for the mass-critical case.

Unlike the energy- and mass-critical problems, for any other  $s_c\notin \{0,1\}$, there are no conserved quantities that control the growth in time of the   $\dot H^{s_c}(\mathbb{R} ^d)$ norm of the solutions.  The first work dealing with Conjecture \ref{CNLS0}  for  $s_c\notin \{0,1\}$ is attributed to  Kenig and Merle \cite{KenigMerle2010} at the case $d = 3, s_c = \frac{1}{2}$ by using their concentration-compactness method developed in \cite{KenigMerle2006}. Further results in the intercritical regime ($0 < s_c < 1$)   can be found in~\cite{GaoMiaoYang2019, GaoZhao2019,LMZ,Murphy2014, Murphy2014b, Murphy2015, XieFang2013, Yu2021}.
 
 In the energy-supercritical case ($s_c > 1$), Killip and Visan \cite{KillipVisan2010} were the first to resolve Conjecture \ref{CNLS0} for
\begin{equation}
	 \begin{cases} 
		1 < s_c < \frac{3}{2} & \text{when } d = 5, 6 \\ 
		1 < s_c < \frac{d + 2 - \sqrt{(d - 2)^2 - 16}}{4} & \text{when } d \geq 7.
	\end{cases} \label{E821}
\end{equation}
 Subsequently, Murphy \cite{Murphy2015} addressed the conjecture for radial initial data in the case $d = 3$ and $s_c \in (1, \frac{3}{2})$.
 By developing long-time Strichartz estimates for the energy-supercritical regime, Miao-Murphy-Zheng \cite{MiaoMurphyZheng2014} and Dodson-Miao-Murphy-Zheng \cite{Dodson2017} resolved the Conjecture \ref{CNLS0} for general initial data when $d = 4$ and $1 < s_c \le \frac{3}{2}$, while the radial case with \(\frac{3}{2} < s_c < 2\) was subsequently addressed by Lu and Zheng \cite{LuZheng2017}.
 More recently, Zhao \cite{Zhao2017AMS} and Li-Li \cite{LiLi2022SIAM} resolved the Conjecture \ref{CNLS0} in the case $d \ge 5$,  where for  $d \ge 8$, their results also required $p$ to be an even number. See Table~\ref{table1} for a summary.

 \begin{table}[h]\label{table1}
 	\centering
 	\caption{Results for Conjecture \ref{CNLS0} in the super-critical case: $ s_c>1$}
 	\begin{tabular}{|c|c|}
 		\hline
 		$d=3$ & $1<s_c<\frac{3}{2}$, \textcolor{blue}{radial}, Murphy \cite{Murphy2015}\\
 		\hline 
 		$d=4$ & \thead {  $1<s_c<\frac{3}{2}$, Miao-Murphy-Zheng\cite{MiaoMurphyZheng2014}; $s_c=\frac{3}{2}$, Dodson-Miao-Murphy-Zheng\cite{Dodson2017}; \\  $\frac{3}{2}<s_c<2$, \textcolor{blue}{radial},  Lu-Zheng\cite{LuZheng2017}}\\
 		\hline 
 		$d\ge5$  & \thead {$1<s_c $, and \textcolor{blue}{  $p$ is even or  $s_c< \frac{d + 2 - \sqrt{(d - 2)^2 - 16}}{4} $ when  $d\ge8$}, \\
 			Killip-Visan\cite{KillipVisan2010}, Zhao\cite{Zhao2017AMS}, Li-Li\cite{LiLi2022SIAM}}\\
 		\hline
 	\end{tabular}
 \end{table}

 In this paper, we  further investigate     Conjecture \ref{CNLS0}  in dimensions  $d\ge8$.  As claimed in \cite{KillipVisan2010}, the results should hold in the whole range of   $s_c>1$, prodvided that   $s_c<1+p$ or  $p$ is an even integer.  The latter case was settled by Li--Li \cite{LiLi2022SIAM}. In this work, we complete the remaining case \(s_c < 1 + p\) by refining some of the nonlinear estimates from \cite{KillipVisan2010,LiLi2022SIAM}.

 	Our main result is stated as follows:
 
 \begin{theorem}\label{T1}
 	Let  $d\ge 8$ and  $1<s_c<1+p$.   
 	Assume that $u: I \times \mathbb{R}  ^d\rightarrow \mathbb{C}$ is a  maximal-lifespan solution to (\ref{NLS}) such that \begin{equation}\label{Ebound}
 		\|u\|_{L_t^{\infty }\dot H_x^{ s_c}(I\times \mathbb{R} ^{d})}\le E<+\infty.
 	\end{equation}
 	Then $u$ is global and scatters as $t \to \pm \infty$ in the sense of (\ref{1.4}).  
 \end{theorem}
  \begin{remark}
 The condition \(s_c < 1+p\) is a natural smoothness assumption, equivalent to  $2p^2 - p(d-2) + 4 > 0.$ 
 This inequality is automatically satisfied for \(d \le 7\), while for \(d \ge 8\) it is weaker than the assumption in (\ref{E821}).
 \end{remark}
 
We conclude this subsection by noting that Conjecture \ref{CNLS0} has also been studied extensively for other variants of the defocusing NLS, including exterior domain problems \cite{LiuSongZheng2026}, tori \(\mathbb{T}^d\) \cite{SongZhang2026,YuYue2024},  fourth order NLS \cite{MiaoZheng2016}, and inhomogeneous NLS \cite{LiuXu2026,WangXu2026}.
 
\subsection{Outline of the proof}

To prove Theorem \ref{T1}, we follow the strategy of \cite{KenigMerle2006}. Suppose, for contradiction, that Theorem \ref{T1} fails. Then one can obtain a special class of solutions that are almost periodic modulo the symmetries of the equation; such results were established in \cite[Theorem 1.12]{KillipVisan2010}.
 
\begin{theorem}[Existence of almost periodic solutions]\label{TReduction}
Suppose  Theorem \ref{T1} fails to be true. Then there exists a maximal-lifespan solution $u: I \times \mathbb{R} ^d\to \mathbb{C}$ to (\ref{NLS})   which blows up in both time directions in the sense that
 \begin{equation}
 	\int _0^ {\sup I} \int _{\mathbb{R} ^d}|u(t,x)|^{\frac{d p}{2}}dxdt= 	\int _ {\inf I}^0 \int _{\mathbb{R} ^d}|u(t,x)|^{\frac{d p}{2}}dxdt=\infty .\label{E811}
 \end{equation}
Moreover,  the solution  $u$ is almost periodic modulo symmetries in the sence that  
there exist functions  $x:I\to \mathbb{R} $,  $N: I \to \mathbb{R}^+$   such that  the set
\begin{equation}
	K=\left\{N(t)^{-\frac{2}{p}}u\left(t,x(t)+\frac{x}{N(t)}\right):t\in I\right\}\label{E12216}
\end{equation}
is precompact in  $\dot H^{s_c}_x (\mathbb{R} ^d)$.   
  \end{theorem}
   \begin{remark}
  	As a consequence of (\ref{E12216}), for any  $\eta >0$, we   can choose $N_0 = N_0(\eta)$ to satisfy
  	\begin{equation}
  		\| |\nabla|^{s_c} P_{\leq N_0} u \|_{L_t^\infty L_x^2} \leq \eta. \label{Ecompact}
  	\end{equation}
  \end{remark}
 
 By Theorem \ref{TReduction}, the proof of Theorem \ref{T1} reduces to ruling out almost periodic solutions. Following the rescaling argument in \cite{KillipVisan2010AJM}, it suffices to rule out almost periodic solutions in the following three cases:
 \begin{itemize}
 	\item[Case 1.] \textit{Finite time blow-up}. $|\inf I| < \infty$ or $|\sup I| < \infty$.
 	\item[Case 2.] \textit{Frequency-cascade}. $u$ is a global solution with $N(t) \geq 1$, and there exists a sequence of time $t_n \to \infty$ such that $N(t_n) \to \infty$ as $n \to \infty$.
 	\item[Case 3.] \textit{Soliton-like solution}. $u$ is a global solution with $N(t) \equiv 1$.
 \end{itemize}
 In the finite time blow-up case, one can apply the  no-waste Duhamel formula (Proposition \ref{Pnowaste} below) together with the Strichartz and Hölder inequalities to show that   $u\equiv0$. This contradicts the fact that  $u$ is a blowup solution.       The details can be found in \cite[Theorem 5.1]{KillipVisan2010} and \cite[Proposition 3.1]{LiLi2022SIAM}.

The exclusion of almost periodic solutions in Cases 2 and 3 relies on the following Proposition. 

\begin{proposition}[Negative regularity]\label{Pnegativeregularity}
	Let  $d\ge8$ and assume the critical regularity   $s_c<1+p$.
	Let  $u$ be the  almost periodic solution in Theorem  \ref{TReduction}, with (\ref{Ebound}) holds.  
	 Then there exists   $\varepsilon >0$  such that   $u\in L_t^{\infty }(\mathbb{R} , \dot H^{-\varepsilon }_x (\mathbb{R} ^d))$.    
\end{proposition}
By invoking Proposition \ref{Pnegativeregularity}, the mass conservation law, and the interaction Morawetz estimate, one  can rule out the almost periodic solutions in Cases 2 and 3. For full details, we refer the reader to \cite[Sections 7 and 8]{KillipVisan2010} and \cite[Section 5]{LiLi2022SIAM}.

Thus, the proof of Theorem \ref{T1} is reduced to establishing Proposition \ref{Pnegativeregularity}, whose proof in turn reduces to proving the following lemma.
  
 \begin{lemma}\label{Lnegativeregularity}
 	Let  $u$  be a solution to (\ref{NLS}) that obeys the hypotheses of Proposition \ref{Pnegativeregularity}.   Then there exists  $l\in (2,\frac{2d}{d-2s_c})$  such that   $u\in L_t^{\infty }L_x^{l}(\mathbb{R}  \times \mathbb{R} ^d)$.   
 \end{lemma}
 
 Indeed,   by (\ref{Ebound}) and Lemma \ref{Lnegativeregularity},    $u\in L_t^{\infty }(I,\dot H_x^{s_c}\cap L^l_x)$. Then    by the double Duhamel formula, and the dispersive  estimate, 
we can prove that if  $u\in L_t^{\infty }(I,\dot H_x^{s}\cap L^l_x)$  with  $s\in [0,s_c]$, then   in fact    $u\in L_t^{\infty }(I,\dot H_x^{s-\varepsilon _0})$, where  $\varepsilon_0 $ is a small    positive constant independent os  $s$.  After a finite number of iterations, we deduce that  $u\in L_t^{\infty }(\mathbb{R} ; \dot H_x^{-\varepsilon }(\mathbb{R} ^d))$.    This completes the proof of Proposition \ref{Pnegativeregularity}.  The  full details can be found in   \cite[Proposition 6.1, Proposition 6.4]{KillipVisan2010} and \cite[Proposition 4.3]{LiLi2022SIAM}. 
 
Thus, to prove Theorem \ref{T1}, it suffices to establish Lemma \ref{Lnegativeregularity}. This lemma was proved for \(d \ge 5\) in \cite[Theorem 5.1]{KillipVisan2010} and \cite[Proposition 4.3]{LiLi2022SIAM}, under   the assumption that \(p\) is even or   $s_c< \frac{d + 2 - \sqrt{(d - 2)^2 - 16}}{4} $ when \(d \ge 8\). In this paper,   we improve the upper bound on \(s_c\) to \(s_c<1+p\) by refining the nonlinear estimates in \cite{KillipVisan2010,LiLi2022SIAM} through a suitable choice of parameters in the argument, treating the different cases separately.

The rest of this paper is organized as follows. In Section~\ref{S2}, we present some notation and useful estimates. Sections~\ref{S3} and~\ref{S4} are devoted to the proof of Lemma~\ref{Lnegativeregularity} in the cases \(s_c \ge 4\) and \(1 < s_c < 4\), respectively.

 \section{Preliminaries}\label{S2}
 
 \subsection{Some notation}
 We use the standard notation for mixed Lebesgue space-time norms and Sobolev spaces. We write  $A\lesssim  B$ to denote  $A\le CB$ for some  $C>0$. If  $A\lesssim B$ and  $B\lesssim A$, then we write  $A\sim B$.  We write  $A\ll B$   to denote  $A\le cB$ for some  small $c>0$.  If  $C$ depends upon some additional parameters, we will indicate this with subscripts; for example,  $X\lesssim _u Y$ denotes that  $X\le C_u Y$ for some  $C_u$ depending on  $u$.    We use  $O(Y)$ to denote any quantity  $X$  such that  $|X|\lesssim  Y$.      We write $L_t^q L_x^r$ to denote the Banach space with norm
 \[  \|u\|_{L_t^q L_x^r (\mathbb{R} \times \mathbb{R}^d)} := \left( \int_{\mathbb{R}} \left( \int_{\mathbb{R}^d} |u(t, x)|^r \, dx \right)^{q/r} \, dt \right)^{1/q},\]
 with the usual modifications when $q$ or $r$ are equal to infinity, or when the domain $\mathbb{R} \times \mathbb{R}^d$ is replaced by spacetime slab such as $I \times \mathbb{R}^d$. When $q = r$ we abbreviate $L_t^q L_x^q$ as $L_{t,x}^q$. 
 
 Let $\phi(\xi)$ be a radial bump function supported in the ball $\{\xi \in \mathbb{R}^d : |\xi| \leq 2\}$ and equal to  $1$ on the ball $\{\xi \in \mathbb{R}^d : |\xi| \leq 1\}$. For each number $N > 0$, we define the Fourier multipliers
 \begin{equation}
 	\widehat{P_{\leq N} f}(\xi) :=\phi(\xi/N) \widehat{f}(\xi), \quad \widehat{P_{> N} f}(\xi) :=(1 - \phi(\xi/N)) \widehat{f}(\xi),\notag
 \end{equation}
 \begin{equation}
 	\widehat{P_N f}(\xi) :=\psi(\xi/N) \widehat{f}(\xi) := (\phi(\xi/N) - \phi(2\xi/N)) \widehat{f}(\xi).\notag
 \end{equation}
 We similarly define $P_{< N}$ and $P_{\geq N}$.  For convenience of notation, let $u_N := P_{N}u$, $u_{\leq N} := P_{\leq N}u$, and $u_{>N} :=P_{>N}u$.

 \subsection{Some useful estimates}
The following three estimates will be used frequently in Sections \ref{S3} and \ref{S4}.
\begin{lemma}[Bernstein estimates]
 For \( 1 \leq r \leq q \leq \infty \),
\[
\||\nabla|^{\pm s} P_N f \|_{L^r_x(\mathbb{R}^d)} \sim N^{\pm s} \|P_N f\|_{L^r_x(\mathbb{R}^d)},
\]
\[
\|P_{\leq N} f\|_{L^q_x(\mathbb{R}^d)} \lesssim N^{\frac{d}{r} - \frac{d}{q}} \|P_{\leq N} f\|_{L^r_x(\mathbb{R}^d)},
\]
\[
\|P_N f\|_{L^q_x(\mathbb{R}^d)} \lesssim N^{\frac{d}{r} - \frac{d}{q}} \|P_N f\|_{L^r_x(\mathbb{R}^d)}.
\]
\end{lemma}

\begin{lemma}[Dispersive estimate]
For \( 2 \leq q \leq \infty \) and \( f \in L^{q'}(\mathbb{R}^d) \),
\[
\left\| e^{it\Delta} f \right\|_{L^q(\mathbb{R}^d)} \lesssim |t|^{-\frac{d}{2}(1-\frac{2}{q})} \left\| f \right\|_{L^{q'}(\mathbb{R}^d)}.
\]
\end{lemma}

 \begin{lemma}[Lemma 2.6 in \cite{Kwak-Kwon}]\label{LF2}
 	Let \(\alpha \geq 1\) and \(m \in \mathbb{Z}\). Let \(F : \mathbb{C} \to \mathbb{C}\) be the function \(F(z) := |z|^{\alpha - m} z^m\). Let \(s \in [0, \alpha)\) and \(p, p_1, p_2 \in (1, \infty)\) be exponents satisfying   \( \frac{1}{p} = \frac{\alpha - 1}{p_1} + \frac{1}{p_2}\). Then, for \(u : \mathbb{R}^d \to \mathbb{C}\), we have
 	\begin{equation} 
 		 \|F(u)\|_{H^{s,p}}\lesssim   \|u\|_{L^{p_1}}^{\alpha - 1}  \|u\|_{H^{s,p_2}}.\notag
 	\end{equation}
 \end{lemma}
 
 The next proposition says in contrast with the classical Duhamel formula that there is no scattered wave at the endpoint of the lifespan \( I \) for almost periodic solutions. We refer to  \cite{TaoVisanZhang2007} for more information.
 
 \begin{proposition}[No-waste Duhamel formula]\label{Pnowaste} 
 	Let \( u : I \times \mathbb{R}^d \to \mathbb{C} \) be a maximal-lifespan almost periodic solution to \eqref{NLS}, then for all \( t \in I \),
 	\begin{align}
 		u(t) &= \lim_{\substack{T \nearrow \sup I}} i \int_{t}^{T} e^{i(t - s)\Delta} (|u|^p u)(s) ds \notag \\
 		&= -\lim_{\substack{T \searrow \inf I}} i \int_{T}^{t} e^{i(t - s)\Delta} (|u|^p u)(s) ds \notag
 	\end{align}
 	as a weak limit in \( \dot{H}_x^{s_c}(\mathbb{R}^d) \).
 \end{proposition}

 Finally, we record  a technical result from \cite[Lemma 2.14]{KillipVisan2010AJM}, which will be  used in Section \ref{S3}.
 
 \begin{lemma}[Acausal Gronwall inequality]\label{LGronwall} 
  Let \( \gamma > 0 \), \( 0 < \eta < \frac{1}{2}(1 - 2^{-\gamma}) \) and \( \{b_k\} \in \ell^\infty(\mathbb{Z}^+) \). Suppose \( \{x_k\} \in \ell^\infty(\mathbb{Z}^+) \) is a nonnegative sequence that satisfies
 	\[
 	x_k \leq b_k + \eta \sum_{\ell = 0}^\infty 2^{-\gamma |k - \ell|} x_\ell
 	\]
 	for all \( k \geq 0 \). Then there exists \( r = r(\eta) \in (2^{-\gamma}, 1) \) such that
 	\[
 	x_k \lesssim \sum_{\ell = 0}^k r^{|k - \ell|} b_\ell
 	\]
 	for all \( k \geq 0 \), with \( r \to 2^{-\gamma} \) as \( \eta \to 0 \).
 \end{lemma}

\section{The case  $s_c\ge4$}\label{S3}

This section is devoted to the proof of Lemma \ref{Lnegativeregularity} in the case \(s_c \ge 4\). Under the assumption \(s_c < 1+p\), we necessarily have \(p > 3\).

Let \(\eta > 0\) be a small parameter to be fixed later, and let \(N_0 = N_0(\eta)\) denote the constant given by (\ref{Ecompact}).

We choose \(\varepsilon > 0\) such that
\[
0 < \varepsilon < \min\left\{\frac14,\, \frac{d - 2s_c}{4}\right\},
\]
and define the exponent 
\[
\frac{1}{q} := \frac{1}{2} - \frac{2 - 4\varepsilon}{d}.
\]

Fix \(N \in 2^{\mathbb{Z}}\). We first derive an estimate on  \(\| u_N(0) \|_{L_x^q}\).  By     Proposition \ref{Pnowaste}, dispersive  estimate and Bernstein's inequality, we have 
 \begin{align}
 		\| u_N(0) \|_{L_x^q} &\leq \left\| \int_0^\infty e^{-it\Delta} P_N F(u)(t) \, dt \right\|_{L_x^q} \notag\\
 	&\leq \left\| \int_0^{N^{-2}} e^{-it\Delta} P_N F(u)(t) \, dt \right\|_{L_x^q} + \left\| \int_{N^{-2}}^\infty e^{-it\Delta} P_N F(u)(t) \, dt \right\|_{L_x^q}\notag \\
 	&\lesssim N^{\frac{d}{2} - \frac{d}{q}} \cdot \int_0^{N^{-2}} \| P_N F(u)(t) \|_{L_x^2} \, dt + \int_{N^{-2}}^\infty t^{-\frac{d}{2}\left(1 - \frac{2}{q}\right)} \| P_N F(u)(t) \|_{L_x^{q'}} \, dt \notag \\
 	&\lesssim N^{\left(1 - \frac{2}{q}\right)d - 2} \| P_N F(u) \|_{L_t^\infty L_x^{q'}}, \label{E12231}
 \end{align}
 where  $F(u):=|u|^pu$.

	 To estimate  $ \| P_N F(u) \|_{L_t^\infty L_x^{q'}}$,   we decompose 
	\[
	F(u) = F(u_{\leq N_0}) + \bigl( F(u) - F(u_{\leq N_0}) \bigr).
	\]
Bernstein's inequality, combined with Hölder's inequality,  Sobolev embedding and the   bound (\ref{Ebound}), gives
	\begin{align}
		&\| P_N \bigl( F(u) - F(u_{\leq N_0}) \bigr) \|_{L_t^\infty L_x^{q'}} \notag\\
		&\lesssim N^{d(\frac{1}{q}-\frac{1}{2})+2} \| F(u) - F(u_{\leq N_0}) \|_{L_t^\infty L_x^{\frac{2d}{d + 4}}} \notag \\
		&\lesssim N^{d(\frac{1}{q}-\frac{1}{2})+2} \| u_{\geq N_0} \|_{L_t^\infty L_x^2} \| u \|_{L_t^\infty L_x^{\frac{dp}{2}}}^p\notag\\
		& \lesssim N^{d(\frac{1}{q}-\frac{1}{2})+2} N_0^{-s_c} \||\nabla|^{s_c}u\|_{L_t^{\infty }L_x^{2}}^{p+1}   \notag\\
		&\lesssim N^{d(\frac{1}{q}-\frac{1}{2})+2} N_0^{-s_c}.   \label{E12232}
	\end{align}
 Now we consider $F(u_{\leq N_0})$. We set $v = u_{\leq N_0}$ and write 
\begin{align}
		F(u_{\le N_0}) &= |v|^p v_{\geq 10N} + |v_{\leq 10N}|^p v_{\leq 10N} + \bigl( |v|^p - |v_{\leq 10N}|^p \bigr) v_{\leq 10N}\notag \\
	&:= F_1(v) + F_2(v) + F_3(v).\label{E12233}
\end{align}

The combination of (\ref{E12231})–(\ref{E12233}) now yields
	\begin{equation}
		\|u_N(0)\|_{L_x^{q}}\lesssim  N^{2-8\varepsilon }[N^{4\varepsilon }N_0^{-s_c}+ \|P_N[F_1(v)+F_2(v)+F_3(v)]\|_{L_t^{\infty }L_x^{q'}}].\label{E1223S}
	\end{equation}
		We estimate the three terms \(F_1(v), F_2(v),\) and \(F_3(v)\) separately below.

	\textit{ { (1) Estimate on  $F_1(v)$. }}
   By H\"older's inequality, we have 
	\begin{align}
		\|P_N (|v|^{p}v_{>10N})\|_{L_x^{q'}}&\lesssim \sum _{M>10N} \|\widetilde{P}_M(|v|^p)v_M\|_{L_x^{q'}}\notag\\
		&\lesssim \sum _{M>10N} \|v_M\|_{L_x^{q}} \|\widetilde{P}_M(|v|^{p})\|_{L_x^{\frac{q}{q-2}}}, \label{E12212}
	\end{align}      
where  $\widetilde{P}_M:=P_{M/2<\cdot <M}$.  Noting that  by (\ref{Es})
\begin{equation}
	\frac{d-2 (s_c-2+8 \varepsilon )}{2d}+(p-1)\frac{2}{d p}=\frac{4-8 \varepsilon }{d}=\frac{q-2}{q}, \notag
\end{equation}   
we deduce  from  Bernstein, Lemma \ref{LF2}, Sobolev embedding and (\ref{Ecompact}) that 
	\begin{align}
		\|\widetilde{P}_M(|v|^{p})\|_{L_x^{\frac{q}{q-2}}}&\lesssim M ^{-(2-8\varepsilon )} \||\nabla |^{2-8\varepsilon }|v|^{p}\|_{L_x^{\frac{q}{q-2}}}\notag\\
		&\lesssim  M^{-(2-8\varepsilon )} \||\nabla |^{2-8\varepsilon }v\|_{L_x^{\frac{2d}{d-2(s_c-2+8\varepsilon )}}} \|v\|_{L_x^{\frac{d p}{2}}}^{p-1}\notag\\
		&\lesssim   M^{-(2-8\varepsilon )}  \||\nabla|^{s_c}u_{\le N_0}\|_{L_x^{2}}^p\notag\\
		&\lesssim \eta ^p M^{-(2-8\varepsilon )}. \label{E12211}
	\end{align}
	Substituting (\ref{E12211}) into (\ref{E12212}), we obtain 
	\begin{equation}
		\|F_1(v)\|_{L_t^{\infty }L_x^{q'}}\lesssim  \eta ^p M^{-(2-8\varepsilon )}\sum _{M\ge 10N} \|v_M\|_{L_x^{q}}.\label{E12a1}
	\end{equation}
	
		\textit{ { (2) Estimate on  $F_2(v)$. }}
  By Bernstein, Lemma \ref{LF2} and Sobolev embedding, we have 
	\begin{align}
		\|P_N(|v_{\le 10N}|^p v_{\le 10N})\|_{L_x^{q'}}&\lesssim  N ^{-(s_c-2\varepsilon )} \||\nabla |^{s_c-2\varepsilon }( |v_{\le 10N}|^p v_{\le 10N})\|_{L_x^{q'}}\notag\\
		&\lesssim  N^{-(s_c-2\varepsilon )} \|v_{\le 10N}\|_{L_x^{\frac{d p}{2-2 \varepsilon }}}^{p} \||\nabla |^{s_c-2\varepsilon }v_{\le 10 N}\|_{L_x^{\frac{2d}{d-4 \varepsilon }}} \notag\\
		&\lesssim N^{-(s_c-2\varepsilon )}   \|v_{\le 10N}\|_{L_x^{\frac{d p}{2-2 \varepsilon }}}^{p}  \||\nabla|^{s_c}u_{\le N_0}\|_{L_x^{2}}  .\label{E12213}
	\end{align} 
	Moreover, by Gagliardo-Nirenberg's inequality and  Sobolev embedding 
\begin{align}
		\|v_{\le 10N}\|_{L_x^{\frac{d p}{2-2 \varepsilon }}}&\lesssim  \|v_{\le 10N}\|_{L_x^{\frac{2d}{d-2s_c-4\varepsilon }}}^{\frac{1}{p}} \| |\nabla|^{s_c}v_{\le 10N}\|_{ L_x^{2}}^{1-\frac{1}{p}}\notag\\
		&\lesssim   \||\nabla|^{s_c-2+6 \varepsilon }v_{\le 10N}\|_{L_x^{ q}}^{\frac{1}{p}} \| |\nabla|^{s_c}v_{\le 10N}\|_{ L_x^{2}}^{1-\frac{1}{p}}\notag .
\end{align}
	Substituting the preceding estimate into (\ref{E12213}) and invoking   (\ref{Ecompact}), we deduce
	\begin{equation}
		\|F_2(v)\|_{L_t^{\infty }L_x^{q'}}\lesssim  \eta^p N ^{-(s_c-2\varepsilon )}\sum _{M\le 10N}M ^{s_c-2+6\varepsilon } \|v_M\|_{L_x^{q}}.\label{E12a2}
	\end{equation}
	
		\textit{ { (3) Estimate on  $F_3(v)$. }}  We proceed as we did for  (\ref{E12213}). This gives 
	\begin{align}
		&\|P_N[(|v|^{p}-|v_{\le 10N}|^p)v_{\le 10N}]\|_{L_x^{q'}}\notag\\
		&\lesssim   \|v_{\le 10N}\|_{L_x^{\frac{2d}{d-2s_c-4\varepsilon }}} \|v_{\ge 10N}\|_{L_x^{\frac{2d}{d-4\varepsilon }}} \|v\|_{L_x^{\frac{d p}{2}}}^{p-1}\notag\\
		&\lesssim  \sum _{M\le 10N}M ^{s_c-2+6\varepsilon } \|v_M\|_{L_x^{q}}N ^{-(s_c-2\varepsilon )} \||\nabla|^{s_c-2 \varepsilon }v\|_{L_x^{\frac{2d}{d-4\varepsilon }}} \|v\|_{L_x^{\frac{d p}{2}}}^{p-1}\notag\\
		&\lesssim   \sum _{M\le 10N}M ^{s_c-2+6\varepsilon } \|v_M\|_{L_x^{q}}N ^{-(s_c-2\varepsilon )}   \||\nabla|^{s_c}u_{\le N_0}\|_{L_x^{2}}^p\notag\\
		&\lesssim \eta ^p N ^{-(s_c-2\varepsilon )}\sum _{M\le 10N}M ^{s_c-2+6\varepsilon } \|v_M\|_{L_x^{q}}.\label{E121}
	\end{align}

Collecting the estimates (\ref{E12a1}), (\ref{E12a2}), and (\ref{E121}) into (\ref{E1223S}) yields
	\begin{equation}
		\|u_N(0)\|_{L_x^{q}}\lesssim  N^{2-4\varepsilon }N_0^{-s_c}+ \eta ^{p} \sum _{M\ge 10N}(\frac{N}{M})^{2-8\varepsilon } \|v_M\|_{L_x^{q}}+ \eta ^p \sum _{M\le 10N}(\frac{M}{N})^{s_c-2+6\varepsilon } \|v_M\|_{L_x^{q}}.\notag
	\end{equation}
The time-translation invariance of the equation allows us to replace the initial datum by the solution at arbitrary time, giving
		\begin{equation}
		\|u_N(t)\|_{L_t^{\infty }L_x^{q}}\lesssim  N^{2-4\varepsilon }N_0^{-s_c}+ \eta ^{p} \sum _{M\ge 10N}(\frac{N}{M})^{2-8\varepsilon } \|v_M\|_{L_t^{\infty }L_x^{q}}+ \eta ^p \sum _{M\le 10N}(\frac{M}{N})^{s_c-2+6\varepsilon } \|v_M\|_{L_t^{\infty }L_x^{q}}.\notag
	\end{equation}

	Let  $A_N=N^{s_c-2+4\varepsilon } \|v_N(t)\|_{L_t^{\infty }L_x^{q}}$.   Then 
	\begin{equation}
		A_N\lesssim  (\frac{N}{N_0})^{s_c} +\eta ^p \sum _{M\ge 10N}(\frac{N}{M})^{s_c-4\varepsilon }A_M+\eta ^p \sum _{M\le 10N}(\frac{M}{N})^{2\varepsilon }A_M.\label{E7281}
	\end{equation}
	We apply Lemma \ref{LGronwall} with  $N=2^{-k}N_0/10, x_k=A(N)$, and  $\eta$ chosen sufficiently small.  Note that  $\{x_k\}\in l^\infty (\mathbb{Z} ^+)$ by  Sobolev embedding and (\ref{Ebound}). Using (\ref{E7281}) and Lemma \ref{LGronwall}, we deduce  
	\begin{equation}
		\||\nabla|^{s_c-2+4\varepsilon }v_M\|_{L_t^{\infty }L_x^{q}}\lesssim M ^{\frac{3}{2}\varepsilon },\notag
	\end{equation}
	which further implies 
	\begin{equation}
		\||\nabla|^{s_c-2+3\varepsilon }u_M\|_{L_t^{\infty }L_x^{q}}\lesssim M^{\frac{1}{2}\varepsilon },\qquad M\le N_0. \label{E122141}
	\end{equation}
	Let 
	\begin{equation}
		\frac{1}{l}:=\frac{1}{q}-\frac{s_c-2+3\varepsilon }{d}=\frac{1}{2}-\frac{ s_c-\varepsilon }{d}.\notag
	\end{equation}
	Then  by Bernstein, Sobolev embedding, (\ref{E122141}) and (\ref{Ebound}),  we have 
	\begin{align}
		\|u\|_{L_t^{\infty }L_x^{l}}&\lesssim  \sum _{M\le N_0} \| |\nabla|^{s_c-2+3\varepsilon }u_M\|_{L_t^{\infty }L_x^{q}}+\sum _{M>N_0}M ^{-\varepsilon } \||\nabla|^{s_c}u\|_{L_t^{\infty }L_x^{2}}\notag\\
		&\lesssim  N_0^{\varepsilon /2}+N_0^{-\varepsilon }. \label{E12215}
	\end{align}
	Hence  $u\in L_t^{\infty }L_x^{l},$ with  $2<l<\frac{2d}{d-2s_c}$.   This completes the proof of Lemma \ref{Lnegativeregularity} in the case  $s_c\ge 4$.

\section{The case  $1<s_c<4$}\label{S4}

This section is devoted to the proof of Lemma \ref{Lnegativeregularity} for \(1 < s_c < 4\). Let \(\eta > 0\) be a small parameter to be fixed later, and let \(N_0 = N_0(\eta)\) be the constant from (\ref{Ecompact}). The argument is divided into two subsections, dealing with the cases \(p > 1\) and \(p \le 1\), respectively. 

 \subsection{The case  $p>1$ }

Choose \(\varepsilon > 0\) sufficiently small such that
\[
0 < \varepsilon < \min \left\{ \frac{1}{3}(p+1-s_c),\, \frac{1}{6}s_c,\, \frac{4-s_c}{12} \right\}.
\]
Set
\[
q := \frac{2d}{d - 2\left(\frac{2+s_c}{3} + \varepsilon\right)},
\qquad
s := \frac{2}{3}(s_c - 1) + 2\varepsilon.
\]
A direct computation shows that  $q \in \left(\frac{2d}{d-2},\, \frac{2d}{d-2s_c}\right)$  and $ s < \frac{s_c}{2}.$

By the same arguments as that used to derive (\ref{E1223S}), we obtain 
\begin{align}
		\|u_N(0)\|_{L_x^{q}} 
		&\lesssim   N ^{(1-\frac{2}{q})d-2} [N ^{d(\frac{1}{q}-\frac{1}{2})+2}N_0^{-s_c}+\|P_N[F_1(v)+F_2(v)+F_3(v)]\|_{L_t^{\infty }L_x^{q'}}]\notag\\
		&\lesssim  N ^{\frac{s}{2}+1}N_0 ^{-s_c}+N^s\|P_N[F_1(v)+F_2(v)+F_3(v)]\|_{L_t^{\infty }L_x^{q'}}, \label{E127}
\end{align}
where  $F_1(v), F_2(v), F_3(v)$ are as in (\ref{E12233}).  

	We now estimate the three terms \(F_1(v), F_2(v),\) and \(F_3(v)\) separately below.

	\textit{ { (1) Estimate on  $F_1(v)$. }}  By Bernstein and H\"older's inequality, we have 
\begin{align}
	\|P_N (|v|^{p}v_{>10N})\|_{L_x^{q'}}&\lesssim \sum _{M>10N} \|\widetilde{P}_M(|v|^p)v_M\|_{L_x^{q'}}\notag\\
	&\lesssim \sum _{M>10N} \|v_M\|_{L_x^{q}} \|\widetilde{P}_M(|v|^{p})\|_{L_x^{\frac{q}{q-2}}},\label{E122}
\end{align}      
where  $\widetilde{P}_M:=P_{M/2<\cdot <M}$.   Noting that  by (\ref{Es})
\begin{equation}
	\frac{d-2(s_c-s)}{2d}+(p-1)\frac{2}{d p}=\frac{s+2}{d}=\frac{q-2}{q},\notag
\end{equation}  
we deduce  from  Bernstein,  Lemma \ref{LF2},  Sobolev embedding and (\ref{Ecompact}) that 
\begin{align}
	\|\widetilde{P}_M(|v|^{p})\|_{L_x^{\frac{q}{q-2}}}&\lesssim M ^{-s} \||\nabla |^{s }|v|^{p}\|_{L_x^{\frac{q}{q-2}}}\notag\\
	&\lesssim  M^{-s}\|v\|_{L_x^{\frac{d p}{2}}}^{p-1} \||\nabla |^{s}v\| _{L_x^{\frac{2d}{d-2(s_c-s)}} }\notag\\
	&\lesssim M^{-s} \||\nabla|^{s_c}u_{\le N_0}\|_{L_x^{2}} ^p \notag\\
	&\lesssim \eta ^{p}M ^{-s}. \label{E123}
\end{align}
Substituting (\ref{E123}) into (\ref{E122}), we obtain 
\begin{equation}
	 \|F_1(v)\|_{L_x^{q}}\lesssim  \eta^p \sum _{M>10N}M ^{-s} \|v_M\|_{L_x^{q}}.\label{E124}
\end{equation}

	\textit{ { (2) Estimate on  $F_2(v)$. }} Let  
	 \begin{equation}
	 	\frac{1}{\overline{q}}:=\frac{1}{q}-\frac{s}{d} \quad\text{and}\quad  \frac{1}{\overline{r}}:=\frac{1}{2}-\frac{s_c-2s}{d}. \label{E7231} 
	 \end{equation}
	 Noting that  by (\ref{Es})
\begin{equation}
	\frac{1}{\overline{q}}+\frac{2}{d p}(p-1)+\frac{1}{\overline{r}}=\frac{1}{q}+\frac{s+2}{d}=\frac{1}{q'}, \notag
\end{equation} 
we deduce from Bernstein, Lemma \ref{LF2}, Sobolev embedding and (\ref{Ecompact})  that 
\begin{align}
 \|F_2(v)\|_{L_x^{q'}}&=	\|P_N(|v_{\le 10N}|^p v_{\le 10N})\|_{L_x^{q'}}\lesssim  N ^{- 2s} \||\nabla |^{2s }( |v_{\le 10N}|^p v_{\le 10N})\|_{L_x^{q'}}\notag\\
	&\lesssim  N^{-2s}   \|v_{\le 10N}\|_{L_x^{\overline{q}}}  \|v_{\le 10N}\|_{L_x^{\frac{d}{2}p}}^{p-1} \||\nabla|^{2s}v_{\le 10N}\|_{L_x^{\overline{r}}}\notag \\
	&\lesssim  N^{-2s}  \|v_{\le 10N}\|_{L_x^{\overline{q}}} \||\nabla|^{s_c}u_{\le N_0}\|_{L_x^{2}}^p \notag\\
	&\lesssim \eta^p N^{-2s}  \|v_{\le 10N}\|_{L_x^{\overline{q}}}\notag\\
	&\lesssim  \eta ^p N^{-2s}\sum _{M\le 10N}M^{s} \|v_M\|_{L_x^{q}}.\label{E125}
\end{align}

	\textit{ { (3) Estimate on  $F_3(v)$. }}   Let  $\overline{q}, \overline{r}$ be defined by (\ref{E7231}).  We proceed as we did for (\ref{E125}). This gives   
\begin{align}
	 \|F_3(v)\|_{L_x^{q'}}&=\|P_N[(|v|^{p}-|v_{\le 10N}|^p)v_{\le 10N}]\|_{L_x^{q'}}\notag\\
	&\lesssim   \|v_{\le 10N}\|_{L_x^{\overline{q}}}\|v_{\ge 10N}\|_{L_x^{\overline{r}}} \|v\|_{L_x^{\frac{d p}{2}}}^{p-1}\notag\\
	&\lesssim  \sum _{M\le 10N}M ^{s } \|v_M\|_{L_x^{q}}N ^{ -2s} \||\nabla|^{2s }v_{\ge 10N}\|_{L_x^{\overline{r}}}  \||\nabla|^{s_c}v\|_{L_x^{2}}^{p-1}\notag\\
	&\lesssim  \sum _{M\le 10N}M^{s} \|v_M\|_{L_x^{q}}N^{-2s}  \||\nabla|^{s_c}u_{\le N_0}\|_{L_x^{2}}^p\notag\\
	&\lesssim  \eta^p N ^{-2s}\sum _{M\le 10N}M ^{s} \|v_M\|_{L_x^{q}}.\label{E126}
\end{align}

Collecting the estimates (\ref{E124})–(\ref{E126}) into (\ref{E127})   gives
\begin{equation}
	\|u_N(0)\|_{L_x^{q}}\lesssim  N^{\frac{s}{2}+1 }N_0^{-s_c}+ \eta ^{p} \sum _{M\ge 10N}(\frac{N}{M})^{s} \|v_M\|_{L_x^{q}}+ \eta ^p \sum _{M\le 10N}(\frac{M}{N})^{s } \|v_M\|_{L_x^{q}}.\notag
\end{equation}
By time-translation symmetry, we also have 
\begin{equation}
	\|u_N(t)\|_{L_t^{\infty }L_x^{q}}\lesssim  N^{\frac{s}{2}+1 }N_0^{-s_c}+ \eta ^{p} \sum _{M\ge 10N}(\frac{N}{M})^{s} \|v_M\|_{L_t^{\infty }L_x^{q}}+ \eta ^p \sum _{M\le 10N}(\frac{M}{N})^{s } \|v_M\|_{L_t^{\infty }L_x^{q}}.\notag
\end{equation}

Let  $A_N=N^{s_c-\frac{s}{2}-1} \|v_N(t)\|_{L_t^{\infty }L_x^{q}}$.  Recal that     $s=\frac{2}{3}(s_c-1)+2 \varepsilon $.   Then 
\begin{equation}
	A_N\lesssim  (\frac{N}{N_0})^{s_c} +\eta^p \sum _{M\ge 10N}(\frac{N}{M})^{ \frac{4}{3}(s_c-1)+\varepsilon  }A_M+ \eta^p \sum _{M\le 10N}(\frac{M}{N})^{3\varepsilon }A_M.\notag
\end{equation}
By the same argument as in (\ref{E122141}), we choose \(\eta > 0\) sufficiently small so that we can apply Lemma \ref{LGronwall} to deduce that
\begin{equation}
	\||\nabla|^{s_c-\frac{s}{2}-1}v_M\|_{L_t^{\infty }L_x^{q}}\lesssim M ^{2\varepsilon },\notag
\end{equation}
which further  implies 
\begin{equation}
	\||\nabla|^{s_c-\frac{s}{2}-1+\varepsilon }u_M\|_{L_t^{\infty }L_x^{q}}\lesssim M^{\varepsilon },\qquad M\le N_0. \label{E12214}
\end{equation}
Let 
\begin{equation}
	\frac{1}{l}:=\frac{1}{q}-\frac{1}{d}(s_c-\frac{s}{2}-1-\varepsilon )=\frac{1}{2}-\frac{s_c-\varepsilon }{d}.\notag
\end{equation}
We then proceed as  in  (\ref{E12215}). This gives 
\begin{equation}
	\|u\|_{L_t^{\infty }L_x^{l}}\lesssim  N_0^{\varepsilon }+N_0^{-\varepsilon }.\notag
\end{equation}
Hence  $u\in L_t^{\infty }L_x^{l}$ with  $2<l<\frac{2d}{d-2s_c}$.     This completes the proof of Lemma \ref{Lnegativeregularity} in the case  $1<s_c< 4, p>1$.

\subsection{The case  $0<p<1$}
 Let  $\varepsilon >0$  sufficiently small  such that 
\begin{equation}
	0<\varepsilon <\min \left\{ p,\ \frac{1}{2}-\frac{(1-p)s_c}{4} \right\}\notag
\end{equation} 
and 
\begin{equation}
	s_1:=ps_c+2-3\varepsilon ,\qquad  s_2:=ps_c-\varepsilon ,\qquad q:=\frac{2d}{d-2\varepsilon }.\notag
\end{equation}
By the same arguments as that used to derive (\ref{E1223S}), we obtain 
\begin{align}
		\|u_N(0)\|_{L_x^{q}} 
	&\lesssim   N ^{(1-\frac{2}{q})d-2} [N ^{d(\frac{1}{q}-\frac{1}{2})+2}N_0^{-s_c}+\|P_N[F_1(v)+F_2(v)+F_3(v)]\|_{L_t^{\infty }L_x^{q'}}]\notag\\
	&\lesssim  N ^{\varepsilon  } N_0 ^{-s_c}+N ^{2(\varepsilon -1)} \|P_N[F_1(v)+F_2(v)+F_3(v)]\|_{L_t^{\infty }L_x^{q'}},\label{E128}
\end{align}
 where  $F_1(v), F_2(v), F_3(v)$ are as in (\ref{E12233}).  

	We now estimate the three terms \(F_1(v), F_2(v),\) and \(F_3(v)\) separately below.

	\textit{ { (1) Estimate on  $F_1(v)$. }} By Bernstein, H\"older, Sobolev embedding  and (\ref{Ecompact}), we have  
\begin{align}
	\|P_N (|v|^{p}v_{>10N})\|_{L_x^{q'}}&\lesssim  N ^{2-2\varepsilon } \|P_N[|v|^{p}v_{\ge 10N}]\|_{L_x^{\frac{2d}{d+4-2\varepsilon }}}\notag\\
	&\lesssim  N ^{2-2\varepsilon }  \|v_{\ge 10N }\|_{L_x^{q}} \|v\|_{L_x^{\frac{d p}{2}}}^p\notag\\
	&\lesssim   N ^{2-2\varepsilon }\sum _{M\ge 10N} \|v_M\|_{L_x^{q}}   \||\nabla|^{s_c}u_{\le N_0}\|_{L_x^{2}}^p\notag\\
	&\lesssim    \eta^p N ^{2-2\varepsilon }\sum _{M\ge 10N} \|v_M\|_{L_x^{q}}. \label{E129}
\end{align}

	\textit{ { (2) Estimate on  $F_2(v)$. }}  In view of the inequality \(s_2 < ps_c < p(p+1) < p+1\), we apply Bernstein's inequality, Lemma \ref{LF2}, the Sobolev embedding, and (\ref{Ecompact}) to obtain
\begin{align}
	\|P_N(|v_{\le 10N}|^p v_{\le 10N})\|_{L_x^{q'}}&\lesssim   N ^{-s_2}\||\nabla |^{ s_2 }[ |v_{\le 10N}|^p v_{\le 10N}]\|_{L_x^{q'}}\notag\\
	&\lesssim  N^{ -s_2} \|v_{\le 10N}\|_{L_x^{\frac{d p}{2  }}}^{p} \||\nabla |^{s_2 }v_{\le 10 N}\|_{L_x^{\frac{2d}{d-4 +2\varepsilon }}}\notag\\
	&\lesssim  N^{ -s_2}  \||\nabla|^{s_c}u_{\le N_0}\|_{L_x^{2}}^p \||\nabla|^{s_1}v_{\le 10N}\|_{L_x^{q}}\notag\\
	&\lesssim  \eta ^p N ^{-s_2}\sum _{M\le 10N}M ^{s_1} \|v_M\|_{L_x^{q}}. \label{E1210}
\end{align}

	\textit{ { (3) Estimate on  $F_3(v)$. }}  Let 
\begin{equation}
	\frac{1}{\overline{q}}:=\frac{1}{q}-\frac{s_1}{d}=\frac{d-2 ps_c-4+4\varepsilon }{2d}\quad\text{and}\quad \frac{1}{\overline{r}}:=\frac{1}{2}-\frac{1}{d}(s_c-\frac{s_2}{p})=\frac{1}{2}-\frac{\varepsilon }{p d}.\notag
\end{equation} 
Since \(s_c < p+1\), \(p < 1\), and \(d \ge 8\), we have
\[
d - 2ps_c - 4 > d - 2p(p+1) - 4 \ge d - 8 \ge 0,
\]
so \(\overline{q} > 0\). Then, by Hölder's inequality, (\ref{Es}), Bernstein's inequality, the Sobolev embedding, and (\ref{Ecompact}), we obtain
\begin{align}
	&\|P_N[(|v|^{p}-|v_{\le 10N}|^p)v_{\le 10N}]\|_{L_x^{q'}}\notag\\
	&\lesssim   \|v_{\le 10N}\|_{L_x^{\overline{q}}}  \|v_{\ge 10N}\|_{L_x^{\overline{r}}}^p \notag\\
	&\lesssim   \||\nabla|^{s_1}v_{\le 10 N}\|_{L_x^{q}} \||\nabla|^{\frac{s_2}{p}}v_{\ge 10N}\|_{L_x^{\overline{r}}}^{p}N^{-s_2}\notag\\
	&\lesssim  N ^{-s_2}\sum _{M\le 10N}M^{s_1} \|v_M\|_{L_x^{q}}  \||\nabla|^{s_c}u\|_{L_x^{2}}^p \notag\\
	&\lesssim  \eta^p  N ^{-s_2}\sum _{M\le 10N}M^{s_1} \|v_M\|_{L_x^{q}} . \label{E1211}
\end{align} 

Substituting (\ref{E129}), (\ref{E1210}) and (\ref{E1211}) into (\ref{E128}), we obtain 
\begin{equation}
	\|u_N(0)\|_{L_x^{q}}\lesssim  N ^{\varepsilon  } N_0 ^{-s_c}+  \eta^p \sum _{M\ge 10N} \|v_M\|_{L_x^{q}} +\eta^p N ^{-s_1}\sum _{M\le 10N}M^{s_1} \|v_M\|_{L_x^{q}} .\notag
\end{equation}
By time-translation symmetry, we also have 
\begin{equation}
	\|u_N(t)\|_{L_t^{\infty }L_x^{q}}\lesssim  N ^{\varepsilon  } N_0 ^{-s_c}+  \eta^p \sum _{M\ge 10N} \|v_M\|_{L_t^{\infty }L_x^{q}} +\eta^p N ^{-s_1}\sum _{M\le 10N}M^{s_1} \|v_M\|_{L_t^{\infty }L_x^{q}} .\notag
\end{equation}

Let  $A_N=N^{s_c-\varepsilon } \|v_M\|_{L_t^{\infty }L_x^{q}}$. Then 
\begin{equation}
	A_N\lesssim  (\frac{N}{N_0})^{s_c}+ \eta^p \sum _{M\ge 10N} (\frac{N}{M})^{s_c-\varepsilon }A_M+\eta^p  \sum _{M\le 10N} (\frac{M}{N})^{(p-1)s_c+2-2 \varepsilon }A_M.\notag
\end{equation} 
By the same argument as in (\ref{E122141}), we choose \(\eta > 0\) sufficiently small so that we can apply Lemma \ref{LGronwall} to deduce that
\begin{equation}
	\||\nabla|^{s_c-\varepsilon }v_M\|_{L_t^{\infty }L_x^{q}}\lesssim M ^{ (p-1)s_c+2-3\varepsilon  },\notag
\end{equation}
which further  implies 
\begin{equation}
	\||\nabla|^{s_c-2\varepsilon }u_M\|_{L_t^{\infty }L_x^{q}}\lesssim M ^{ (p-1)s_c+2-4\varepsilon  }, \qquad M\le N_0. \label{E12214}
\end{equation}
Let 
\begin{equation}
	\frac{1}{l}:=\frac{1}{q}-\frac{s_c-2\varepsilon }{d}= \frac{1}{2}-\frac{s_c-\varepsilon }{d}.\notag
\end{equation}
We then proceed as  in  (\ref{E12215}). This gives 
\begin{equation}
	\|u\|_{L_t^{\infty }L_x^{l}}\lesssim  N_0 ^{(p-1)s_c+2-4\varepsilon }+N_0^{-\varepsilon }.\notag
\end{equation}
Hence  $u\in L_t^{\infty }L_x^{l}$ with  $2<l<\frac{2d}{d-2s_c}$.    This completes the proof of Lemma \ref{Lnegativeregularity} in the case  $1<s_c< 4, 0<p< 1$.  

\vspace{1cm}

{\bf Acknowledgements:}   
X.~Liu was supported by NSFC Grant 12501314.


\begin{thebibliography}{99}
 
\bibitem{Bourgain1999} J. Bourgain, \emph{Global well-posedness of defocusing critical nonlinear Schr\"odinger equation in the radial case,} J. Amer. Math. Soc. \textbf{12} (1999), 145-171.
		
 
		
\bibitem{CazenaveWeissler1990NA} T. Cazenave, F. B. Weissler,   \emph{The Cauchy problem for the critical nonlinear Schr\"odinger equation in $H^s$},   Nonlinear Anal. \textbf{14} (1990), no.10, 807–836.
	 
	 
\bibitem{Colliander2008} J. Colliander, M. Keel, G. Staffilani, H. Takaoka, T. Tao, \emph{Global well-posedness and scattering for the energy-critical nonlinear Schr\"odinger equation in $\mathbb{R}^3$,} Ann. of Math. \textbf{167} (2008), 767--865.
	
		
\bibitem{Dodson2012} B. Dodson, \emph{Global well-posedness and scattering for the defocusing, $L^2$-critical nonlinear Schr\"odinger equation when $d\ge3$,} J. Amer. Math. Soc. \textbf{25} (2012), 429--463.
		
\bibitem{Dodson2015} B. Dodson, \emph{Global well-posedness and scattering for the mass critical nonlinear Schr\"odinger equation with mass below the mass of the ground state,} Adv. Math. \textbf{285} (2015), 1589--1618.
		
\bibitem{Dodson2016a} B. Dodson, \emph{Global well-posedness and scattering for the defocusing $L^2$-critical nonlinear Schr\"odinger equation when $d=2$,} Duke Math. J. \textbf{165} (2016), 3435--3516.
		
\bibitem{Dodson2016b} B. Dodson, \emph{Global well-posedness and scattering for the defocusing $L^2$-critical nonlinear Schr\"odinger equation when $d=1$,} Am. J. Math. \textbf{138} (2016), 531--569.
		
	 
\bibitem{Dodson2017} B. Dodson, C. X. Miao, J. Murphy, J. Zheng, \emph{The defocusing quintic NLS in four space dimensions,} Ann. Inst. H. Poincar\'e Anal. Non Lin\'eaire \textbf{34} (2017), 759--787.
		
		
\bibitem{GaoMiaoYang2019} C. Gao, C. Miao, J. Yang, \emph{The Intercritical Defocusing Nonlinear Schr\"odinger Equations with Radial Initial Data in Dimensions Four and Higher,} Anal. Theory Appl. \textbf{35} (2019), 205--234.
		
\bibitem{GaoZhao2019} C. Gao, Z. Zhao, \emph{On scattering for the defocusing high dimensional inter-critical NLS,} J. Differential Equations \textbf{267} (2019), 6198--6215.
			 
 
\bibitem{Kwak-Kwon} B. Kwak, S. Kwon, \emph{Critical local well-posedness of the nonlinear Schr\"odinger equation on the torus}, Ann. Inst. H. Poincar\'e Anal. Non Lin\'eaire \textbf{43} (2024), 155--201. 
	 
\bibitem{KenigMerle2006} C. E. Kenig, F. Merle, \emph{Global well-posedness, scattering and blow-up for the energy-critical, focusing, non-linear Schr\"odinger equation in the radial case,} Invent. Math. \textbf{166} (2006), no. 3, 645-675.
		
\bibitem{KenigMerle2010} C. E. Kenig, F. Merle, \emph{Scattering for $\dot{H}^{1/2}$ bounded solutions to the cubic, defocusing NLS in 3 dimensions,} Trans. Am. Math. Soc. \textbf{362} (2010), 1937--1962.
		
				
\bibitem{KillipTaoVisan2009} R. Killip, T. Tao, M. Visan, \emph{The cubic nonlinear Schr\"odinger equation in two dimensions with radial data,} J. Eur. Math. Soc. \textbf{11} (2009), 1203--1258.

\bibitem{KillipVisan2010AJM} R. Killip, M. Visan, \emph{The focusing energy-critical nonlinear Schr\"odinger equation in dimensions five and higher}, Am. J. Math. \textbf{132} (2010), 361--424.
		
\bibitem{KillipVisan2010} R. Killip, M. Visan, \emph{Energy-supercritical NLS: Critical $\dot{H}^s$-bounds imply scattering,}
Comm. Partial Differential Equations \textbf{35} (2010), 945--987.
		
 
		
\bibitem{KillipVisanZhang2008} R. Killip, M. Visan, X. Zhang, \emph{The mass-critical nonlinear Schr\"odinger equation with radial data in dimensions three and higher,} Anal. PDE \textbf{1} (2008), 229--266.
		
\bibitem{LiLi2022SIAM} J. Li,  K. Li,  \emph{The Defocusing Energy-supercritical Nonlinear Schr\"odinger Equation in High Dimensions},  SIAM J. Math. Anal. \textbf{54} (2022), no.3, 	 3253-3274. 
	
		\bibitem{LMZ} X. Liu, C. Miao, J. Zheng, \emph{ Global well-posedness and scattering for the defocusing one dimensional NLS with algebraica nonlinearity,
			 }  submitted. 
			 
			 \bibitem{LiuSongZheng2026} X. Liu, Y. Song, J. Zheng, \emph{Scattering theory for the defocusing 3d NLS in the exterior of a strictly convex obstacle}, Bull. Sci. Math. \textbf{213} (2026), Paper No. 103875.
			 
			 \bibitem{LiuXu2026} X. Liu, C. Xu, \emph{The defocusing energy-supercritical inhomogeneous NLS in four space dimension}, J. Math. Anal. Appl. \textbf{554} (2026), 129968.  
 
\bibitem{LuZheng2017} C. Lu, J. Zheng, \emph{The radial defocusing energy-supercritical NLS in dimension four,} J. Differ. Equ. \textbf{262} (2017), 4390--4414.
		
\bibitem{MiaoMurphyZheng2014} C. Miao, J. Murphy, J. Zheng, \emph{The defocusing energy-supercritical NLS in four space dimensions,} J. Funct. Anal. \textbf{267} (2014), 1662--1724.
				
		\bibitem{MiaoZheng2016} C. Miao, J. Zheng, \emph{Scattering theory for the defocusing fourth-order Schr\"odinger equation}, Nonlinearity \textbf{29} (2016), 692--736.
		
\bibitem{Murphy2014} J. Murphy, \emph{Inter-critical NLS: Critical $\dot{H}^s$-bounds imply scattering,} SIAM J. Math. Anal. \textbf{46} (2014), 939--997.
		
\bibitem{Murphy2014b} J. Murphy, \emph{The defocusing $\dot{H}^{1/2}$-critical NLS in high dimensions,} Discrete Contin. Dyn. Syst. Ser. A \textbf{34} (2014), 733--748.
		
\bibitem{Murphy2015} J. Murphy, \emph{The radial defocusing nonlinear Schr\"odinger equation in three space dimensions,}
Comm. Partial Differential Equations \textbf{40} (2015), 265--308.
		
\bibitem{RyckmanVisan2007} E. Ryckman, M. Visan, \emph{Global well-posedness and scattering for the defocusing energy critical nonlinear Schr\"odinger equation in $\mathbb{R}^{1+4}$,} Am. J. Math. \textbf{129} (2007), 1--60.

\bibitem{SongZhang2026} Y. Song, R. Zhang, \emph{Global well-posedness for the defocusing cubic nonlinear Schr\"odinger equation on $\mathbb{T}^3$}, Forum Math. \textbf{38} (2026), 809--839.
				
\bibitem{TaoVisanZhang2007} T. Tao,  M. Visan, X. Zhang,  \emph{Global well-posedness and scattering for the mass-critical nonlinear Schr\"odinger equation for radial data in high dimensions,} Duke Math. J. \textbf{140} (2007),  165--202. 
				
\bibitem{Visan2007} M. Visan, \emph{The defocusing energy-critical nonlinear Schr\"odinger equation in higher dimensions,}Duke Math. J. \textbf{138} (2007), 281--374.
		
\bibitem{Visan2012} M. Visan, \emph{Global well-posedness and scattering for the defocusing cubic nonlinear Schr\"odinger equation in four dimensions,} Int. Math. Res. Not. IMRN \textbf{2012} (2012), 1037--1067.
				
		\bibitem{WangXu2026} Y. Wang, C. Xu, \emph{Defocusing $\dot{H}^{\frac{1}{2}}$-critical inhomogeneous nonlinear Schr\"odinger equations}, J. Math. Anal. Appl. \textbf{554} (2026), 129968.
		
\bibitem{XieFang2013} J. Xie, D. Fang, \emph{Global well-posedness and scattering for the defocusing $\dot{H}^s$-critical NLS,} Chin. Ann. Math. \textbf{34B} (2013), 801--842.
		
\bibitem{Yu2021} X. Yu, \emph{Global well-posedness and scattering for the defocusing $\dot{H}^{1/2}$-critical nonlinear Schr\"odinger equation in $\mathbb{R}^2$,} Anal. PDE \textbf{14} (2021), 1037--1067.
		
		\bibitem{YuYue2024} X. Yu, H. Yue, \emph{On the global well-posedness for the periodic quintic nonlinear Schr\"odinger equation}, SIAM J. Math. Anal. \textbf{56} (2024), 1851--1902.
	
\bibitem{Zhao2017AMS} T. Zhao,  \emph{The Defocusing Energy-supercritical NLS in Higher Dimensions},  Acta Mathematica Sinica, English Series,  \textbf{33} (2017), 911--925. 
\end{thebibliography}
\end{document}